\documentclass[12pt]{article}
\usepackage{amsfonts,amssymb,eucal,amsmath}
\usepackage{graphicx}
\title{\bf The Cost of Two-dimensional Rearrangement}
\author{\bf {\small\sc Hailun Zhou}\\
{\small\sc Department of Mathematics, Fudan University,}\\
{\small\sc Shanghai, 200433, China}}
\date{}

\begin{document}
\maketitle

\vspace{0.2in}

\noindent {\sc Introduction} Consider a two dimensional torus 
$\mathbb{T}^2 = \mathbb{R}^2/\mathbb{Z}^2$, with coordinates 
$x = (x_1,x_2) \in [0,1) \times [0,1)$. Let  $A = \{(x_1,x_2)| \ 0 
\leqslant x_2 < 1/2\} \subset \mathbb{T}^2$ be a subset, 
an diffeomorphism $\Phi: \mathbb{T}^2 \to \mathbb{T}^2$ is 
called a rearrangment of $A$.

\vspace{0.1in}

\begin{figure}[h]
\begin{center}
\includegraphics[width=5in]{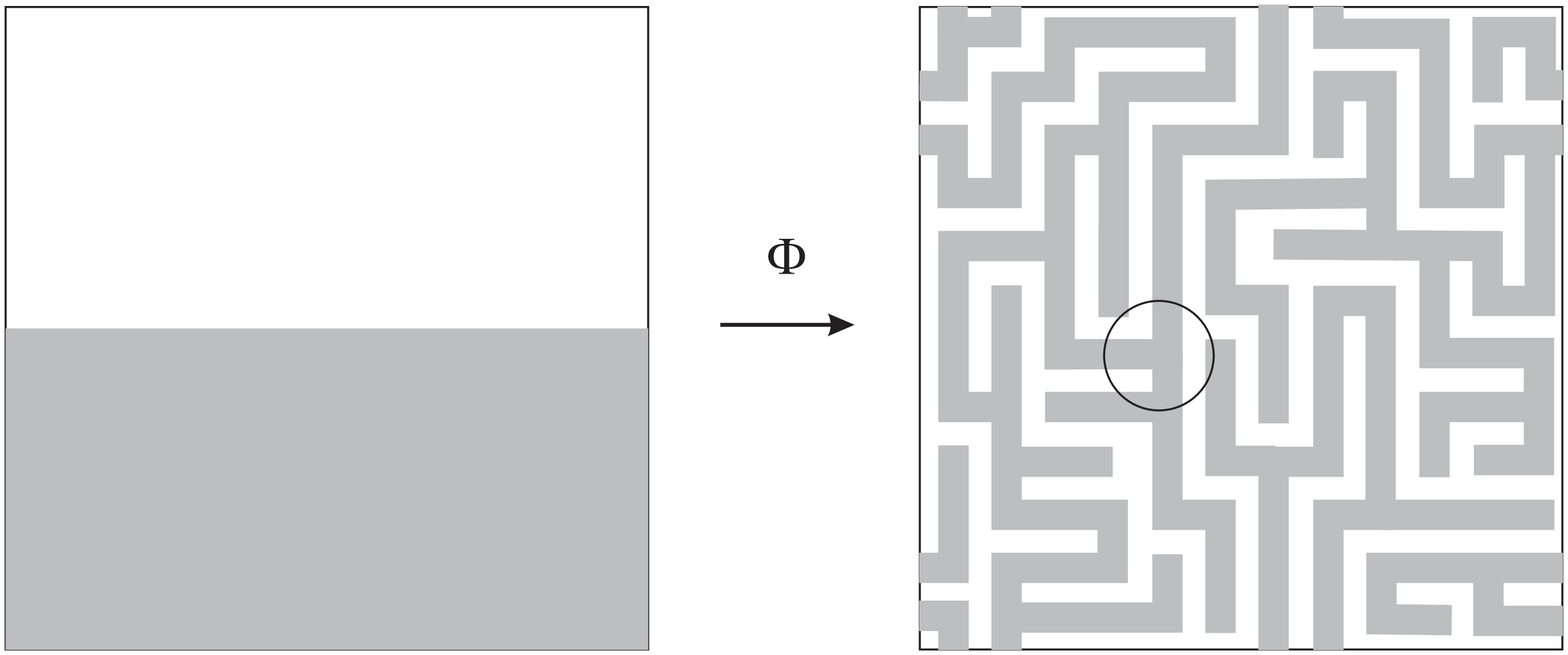}
\end{center}
\end{figure}

\vspace{-0.2in}

\begin{center}
Figure 1
\end{center}

We say that $\Phi$ {\em mixes the set $A$ up to scale} $\varepsilon$ 
if the following holds: there is a fixed real number $\kappa \in (0,1/2)$, 
for any ball $B_\varepsilon(x)$ centered at a point $x \in \mathbb{T}^2$ 
with radius $\varepsilon$, we have
\begin{equation}
\kappa \ \mbox{Area} (B_\varepsilon(x)) \leqslant \mbox{Area} (B_\varepsilon(x) 
\cap \Phi(A)) \leqslant (1-\kappa) \mbox{Area} (B_\varepsilon(x))
\end{equation}

For the rearrangement $\Phi$, let
$$
e(\Phi) = \frac{1}{2} [(\Phi^1_{x_1})^2 + (\Phi^1_{x_2})^2 
+ (\Phi^2_{x_1})^2 + (\Phi^2_{x_2})^2]
$$ 
be the {\em energy density}, and we define the cost of the rearrangement 
as the energy of $\Phi$: 
$$
E(\Phi) = \int_{\mathbb{T}^2} e(\Phi) d \sigma.
$$

The main result in this paper is the following theorem.

{\bf Theorem}. Let  $\Phi(x): \mathbb{T}^2 \to \mathbb{T}^2$ is a 
diffeomorphism and it {\em mixes the set $A$ up to scale} $\varepsilon$. 
If $\Phi$ satisfies
\begin{equation}
0 < \kappa' \leqslant \left| \mbox{det}[\nabla_x \Phi] \right|
\end{equation}
Then there exists a constant $C$ which depends on $\kappa, \kappa'$ only, 
such that
\begin{equation}
E(\Phi) \geqslant \frac{C}{ \varepsilon^2}.
\end{equation}

If $F: [0,1] \times \mathbb{T}^2 \to \mathbb{R}^2 $ is a time 
dependent smooth vector field on $ \mathbb{T}^2 $, and $ \Phi_{t}: \mathbb{T}^2 
\mapsto \mathbb{T}^2$ the flow associated with the vector field 
$F$, i.e., $\Phi_{t}$ is the solution of the following initial value problem
\begin{equation}
\left \{ \begin{array}{ll}
\dot{\Phi}_t(x)=F(t,\Phi_t(x))\\
\Phi_0(x)=x
\end{array}
\right. .
\end{equation}
Let $\Phi_1(\cdot) = {\Phi}(\cdot)$ be the value of the flow at time $t = 1$. 

In [Br], Bressan made the following conjecture:

{\bf Conjecture}. If the flow $\Phi_t$ generated by smooth vector field $F$ is 
{\em nearly incompressible}, i.e., for some constant 
$\kappa' > 0$, we have
\begin{equation}
\kappa' \mbox{Area}(\Omega) \leqslant \mbox{Area} (\Phi_t(\Omega)) 
\leqslant \frac{1}{\kappa'} \mbox{Area}(\Omega),
\end{equation}
for any measurable set $\Omega \subset \mathbb{T}^2$ and $t \in [0,1]$, 
and $\Phi$ {\em mixes the set $A$ up to scale $\varepsilon$}, then there 
is a constant $C$ depends on $\kappa$ and $\kappa'$ only, such that 
\begin{equation}
\int_0^1 \int_{{\mathbb T}^2} | \nabla_x F | d \sigma dt 
\geqslant C | \log \varepsilon |.
\end{equation}

As a corollary, we will prove the following corollary which is in the same manner
as Bressan's conjecture.

{\bf Corollary}. Let $F = F(t, x)$ be a smooth vector field on $\mathbb{T}^2$, and 
assume that the associated flow $\Phi_t(x)$ satisfies
\begin{equation}
0 < \kappa' \leqslant | \mbox{det}[\nabla_x \Phi_t] |,
\end{equation}
and $\Phi = \Phi_1$ {\em mixes the set $A$ up to scale $\varepsilon$}.
Then there exists constant $C$ depends on $\kappa, \kappa'$ only, such that 
\begin{equation}
\frac{C}{\varepsilon^2} \leqslant \int_0^1 \int_{\mathbb{T}^2}
e^{\sqrt6 | \nabla_x F |} d \sigma dt
\end{equation}

\vspace{0.2in}

\begin{center}
{\sc Proof of the Main Theorem and Corollary}
\end{center}

We first prove the theorem.

{\bf Theorem}. Let  $\Phi(x): \mathbb{T}^2 \to \mathbb{T}^2$ is a 
diffeomorphism and it {\em mixes the set $A$ up to scale} $\varepsilon$. 
If $\Phi$ satisfies
\begin{equation}
0 < \kappa' \leqslant \left| \mbox{det}[\nabla_x \Phi] \right|
\end{equation}
Then there exists a constant $C$ which depends on $\kappa, \kappa'$ only, 
such that
\begin{equation}
E(\Phi) \geqslant \frac{C}{ \varepsilon^2}.
\end{equation}

{\bf Proof}. On the set $A \subset \mathbb{T}^2$, we have
$$
\int_A e(\Phi) d \sigma \geqslant \frac{1}{2} \int^{1/2}_0
\left(  \int^1_0 [(\Phi^1_{x_1})^2 + (\Phi^2_{x_1})^2] dx_1 \right) dx_2,
$$
and by H\"{o}lder inequality we have
$$
\int^1_0 [(\Phi^1_{x_1})^2+(\Phi^2_{x_1})^2] dx_1 \geqslant
\left( \int^1_0 \left[ (\Phi^1_{x_1})^2+(\Phi^2_{x_1})^2 \right]^{1/2} dx_1 \right)^2
$$
If we fix $x_2 = s$, then we can view $\Phi(x_1,s)$ as a curve $C_s \subset
\mathbb{T}^2$. Thus,
$$
l(s) = \int^1_0 \left[ (\Phi^1_{x_1}(x_1, s))^2 + (\Phi^2_{x_1}(x_1, s))^2
\right]^{1/2} dx_1
$$
is the length of the curve $C_s$. Then 
\begin{equation*}
\begin{split}
\int_A e(\Phi) & \geqslant \frac{1}{2} \int^{1/2}_0 l^2(s) ds
= \frac{1}{2} \int^{1/4}_0 \left[ l^2(s) + l^2 \left( \frac{1}{2} - s \right)
\right] ds \\
& \geqslant \frac{1}{4} \int^{1/4}_0 \left[ l(s) + l \left( \frac{1}{2} - s \right) 
\right ]^2 ds.
\end{split}
\end{equation*}

Let $A_s = \{ (x_1,x_2) | \ s \leqslant x_2 \leqslant 1/2 - s \} \subset 
\mathbb{T}^2$, and $l(s) + l(1/2 - s)$ is just the length of boundary 
$ \partial(\Phi(A_s))$. Now we will estimate the length of $\partial(\Phi(A_s))$ 
in terms of $\varepsilon$.

\begin{figure}[h]
\begin{center}
\includegraphics[width=2in]{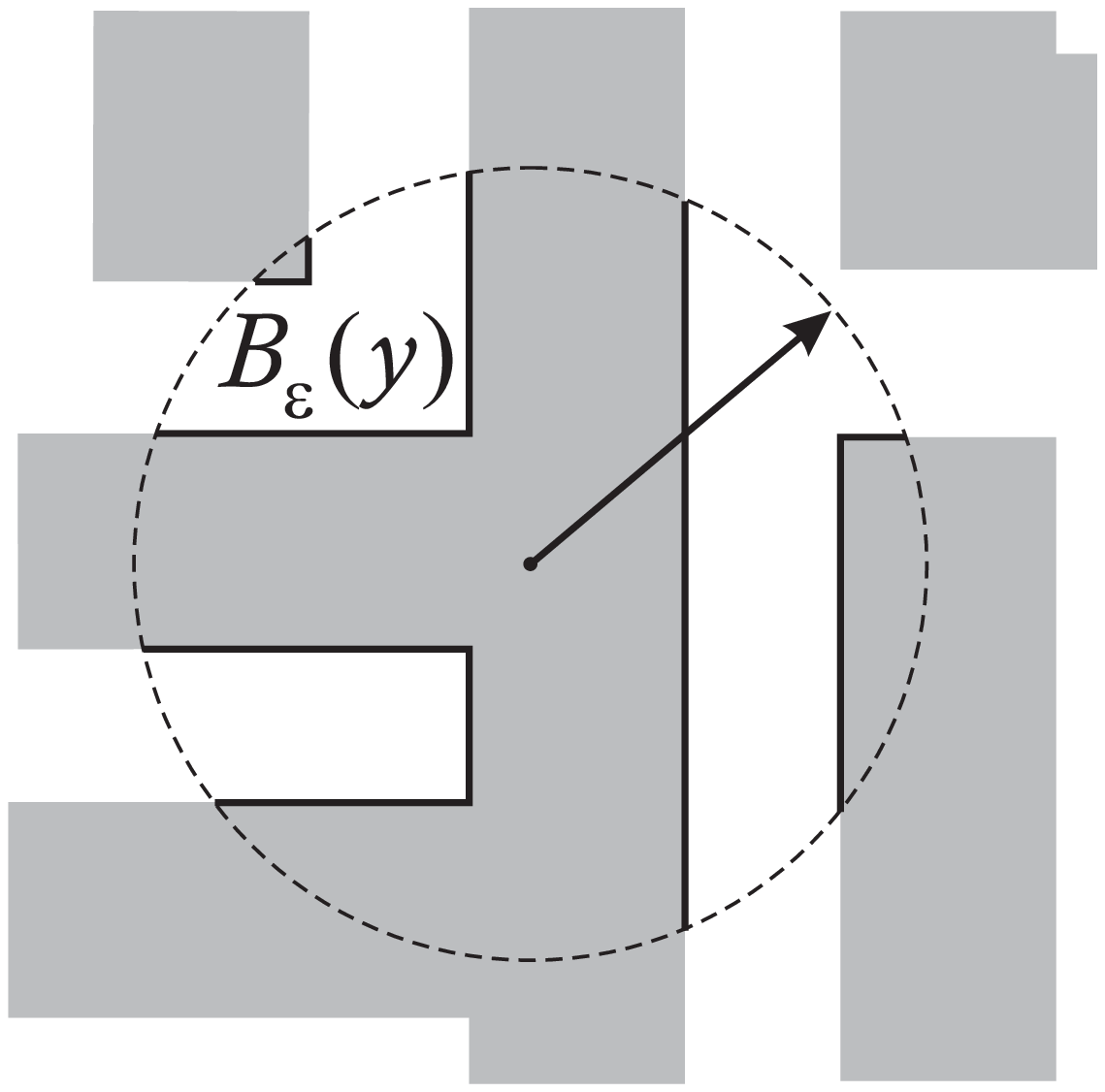}
\end{center}
\end{figure}

For a given point $y \in \Phi (A_s)$, let $B_\varepsilon(y)$ be the
ball centered at $y$ with radius $\varepsilon$, if 
$\kappa \pi \varepsilon^2 > \mbox{Aera} 
(\Phi(A_s) \cap B_{\varepsilon}(y))$, then for any $r > \varepsilon / \sqrt{2}$, 
$\partial B_{r}(y)$ has at least two points intersects with
$\partial (\Phi(A_s))$, and this will imply that
\begin{equation*}
\mbox{length}(\partial (\Phi(A_s)) \cap B_{\varepsilon}(y)) 
\geqslant 2(\varepsilon - \varepsilon / \sqrt{2}) = (2 - \sqrt{2}) \varepsilon.
\end{equation*}
Otherwise, $\partial (B_{r}(y)) \cap \partial(\Phi(A_s)) = \varnothing$. Two curves in
$\partial (\Phi(A_s))$ are both homologically non-trivial, and $y \in B_{r}(y) 
\cap \Phi(A_s)$, so we must have $B_{r}(y) \subset \Phi(A_s)$. This 
contradicts to the assumption $\kappa \pi \varepsilon^2 > \mbox{Aera} (
\Phi(A_s) \cap B_{\varepsilon}(y))$.

\vspace{-0.2in}

\begin{center}
Figure 2
\end{center}

If $\kappa \pi \varepsilon^2 \leqslant \mbox{Area} (\Phi(A_s) \cap B_{\varepsilon}(y))$,
by highly mixing condition on $\Phi$, we have 
\begin{equation*}
\kappa \pi \varepsilon^2 \leqslant \mbox{Area}(\Phi(A_s) \cap B_{\varepsilon}(y)) 
\leqslant \mbox{Area}(\Phi(A) \cap B_{\varepsilon}(y)) 
\leqslant (1 - \kappa) \pi \varepsilon^2.
\end{equation*}
The minimal curve which separates two regions of areas 
$\mbox{Area}(\Phi(A_s) \cap B_{\varepsilon}(y))$ and $\pi \varepsilon^2 -
\mbox{Area}(\Phi(A_s) \cap B_{\varepsilon}(y))$ is a circular arc perpendicular to 
$\partial B_{\varepsilon}(y)$, so there is a constant $m'_{\kappa}$ depends on 
$\kappa$ only, such that 
\begin{equation*}
\mbox{length}(\partial (\Phi(A_s)) \cap B_{\varepsilon}(y)) \geqslant m'_{\kappa} \varepsilon.
\end{equation*} 

Combine these two cases together, let $m_{\kappa} = \min 
\{ m_{\kappa}', (2- \sqrt{2}) \}$, we get a low bound estimation for the length of 
$\partial(\Phi(A_s)) \cap B_{\varepsilon}(y)$: 
\begin{equation*}
\mbox{length}(\partial(\Phi(A_s) \cap B_{\varepsilon}(y)) \geqslant m_{\kappa} \varepsilon.
\end{equation*}

Now we pack the set $\Phi(A_s)$ by a maximal set of balls 
$\{B_{\varepsilon}(y_i) |\ y_i \in \Phi(A_s) \}$ and any two balls 
in the set are disjoint. Let $n$ be the number of balls in this 
maximal set. We have
$$
l(s) + l \left( \frac{1}{2}-s \right) \geqslant m_{\kappa} n \varepsilon.
$$

On the other hand, we notice that balls $\{B_{2\varepsilon}(y_i)\}$
will cover $\Phi(A_s)$. If not, suppose $y_0 \in \Phi(A_s)$ cannot be 
covered by $\{ B_{2\varepsilon}(y_i) \}$,  that means the distance 
between $y_0$ and all the $B_{\varepsilon}(y_i)$ is larger 
than $\varepsilon$. It contradict the maximality of $\{B_{\varepsilon}(y_i) 
|\ y_i \in \Phi(A_s) \}$. Thus we have
$$
4 \pi \varepsilon^2 n \geqslant \mbox{Aera} (\Phi(A_s)).
$$

By condition that $\kappa' < |\mbox{det} (\nabla_x \Phi) |$, we have 
$$
\mbox{Aera} (\Phi(A_s)) \geqslant \kappa' \mbox{Area}(A_s) = 
\kappa' (1/2 - 2s).
$$
Then
$$
4 \pi \varepsilon^2 n \geqslant \kappa'(1/2-2s).
$$
Hence
$$
l(s)+l \left( \frac{1}{2} -s \right) \geqslant \frac{ \kappa' m_{\kappa} }{ 4 \pi \varepsilon }
\left( \frac{1}{2} - 2s \right),
$$
and then
\begin{equation*}
\int_A e(\Phi) d \sigma \geqslant \frac{1}{4} \int_0^{\frac{1}{4}}
\left[ \frac{ \kappa' m_{\kappa}}{ 4 \pi \varepsilon } \left( \frac{1}{2} - 2s
\right) \right]^2 ds = \frac{1}{48} \left(\frac{\kappa' m_{\kappa}}{8
\pi \varepsilon} \right)^2.
\end{equation*}

When $\Phi$ mixes set $A$ up to scale $\varepsilon$, it mixes 
$\mathbb{T}^2 - A$ as well. Similarly, we have
$$
\int_{\mathbb{T}^2-A} e(\Phi) \geqslant 
\frac{1}{48}\left(\frac{M_{\kappa} \kappa'}{8\pi\varepsilon}\right)^2.
$$
Thus
\begin{equation}
\int_{\mathbb{T}^2} e(\Phi) \geqslant 
\frac{1}{24}\left(\frac{M_{\kappa} \kappa'}{8\pi\varepsilon}\right)^2
=\frac{C}{\varepsilon^2}
\end{equation}
where $C$ depends on $\kappa, \kappa'$ only, and we complete the proof
of the theorem. \ \rule{2mm}{3mm}

Now we turn to the proof of the corollary.

{\bf Proof of the Corollary}. We take differentiation with respect to $x$ 
on both sides of the ordinary differential equation $\dot{\Phi}_t(x)=F(t,\Phi_t(x))$. 
According to chain rule, we get 
$$
\partial_t ( \nabla_x \Phi_t(x)) = \nabla_x F(t, \Phi_t(x))
\nabla_x \Phi_t(x).
$$

Let
$$
|\nabla_x F|=\left( \sum^2_{i, j = 1}( \partial_{x_i} F^j )^2 \right)^{1/2}
$$
be the variation of the vector field $F$. Then we have
\begin{equation*}
\partial_t e(\Phi_t) = \sum^2_{i, j=1}( \partial_{x_i} \Phi_t^j )
( \partial_t (\partial_{x_i} \Phi_t^j )) \leqslant 
\sqrt6 |\nabla_x F (t, \Phi_t(x)) | e(\Phi_t ).
\end{equation*}
Then
$$
\partial_t \log e(\Phi_t) = \frac{ \partial_t e(\Phi_t) }{e(\Phi_t)} 
\leqslant \sqrt6 | \nabla_x F(t, \Phi_t(x))|.
$$
Integrating over 0 to 1 from both side we obtain
$$
\log e(\Phi) \leqslant \sqrt6\int^1_0 | \nabla_x F (t, \Phi_t(x))| dt.
$$
By Jensen's inequality
$$
e(\Phi) \leqslant \exp \left[ \sqrt6\int^1_0 | \nabla_x F(t, \Phi_t(x)) | dt \right]
\leq \int^1_0 e^{\sqrt6 | \nabla_x F (t, \Phi_t(x)) |} dt
$$
Integrating over $\mathbb{T}^2$ we have
\begin{equation}
\begin{split}
\int_{\mathbb{T}^2} e(\Phi) d \sigma & \leqslant \int_0^1 \int_{\mathbb{T}^2} 
e^{ \sqrt6 | \nabla_x F(t, \Phi_t(x)) |} d \sigma dt \\
& \leqslant \int_0^1 \int_{\mathbb{T}^2} e^{ \sqrt6 | \nabla_x F(t, \Phi_t(x)) |} 
|\mbox{det}( \nabla_x \Phi_t )|^{-1} d \sigma dt \\
& \leqslant \frac{1}{\kappa'} \int_0^1 \int_{\mathbb{T}^2} 
e^{ \sqrt6 | \nabla_x F(t, x) |} d \sigma dt.
\end{split}
\end{equation}
This completes the proof of the corollary. \rule{2mm}{3mm}

\vspace{0.2in}

\noindent {\sc Reference}

[Br] A. Bressan, A lemma and a conjecture on the cost of rearrangements,
{\sl Rend. Sem. Mat. Univ. Padova}, {\bf 110}(2003), 97-102.

\end{document}